\documentclass[12pt]{article}
\usepackage{amsmath,amssymb,url}
\usepackage[none]{hyphenat}
\usepackage{graphicx}
\begin{document}\begin{center}{\Large{\bf
Dickman multiple polylogarithms\\
and the Lindemann--Furry letters}\\
David Broadhurst\footnote{
School of Physical Sciences, Open University,
Milton Keynes MK7 6AA, UK.
{\tt David.Broadhurst@open.ac.uk}}
and Stephan Ohlmeyer\footnote{
Harvard Business School,
Boston MA 02163, USA.
{\tt sohlmeyer@amp204.hbs.edu}}\\
30 April 2023}\end{center}{\large

{\bf Abstract:} The Dickman function $\rho(u)$ gives the asymptotic probability 
that a large integer $N$ has no prime divisor exceeding $N^{1/u}$.
We expand it in terms of rapidly computable multiple polylogarithms, as defined
by Goncharov and intensively used for evaluations of Feynman integrals
in quantum field theory. In parallel, we solve Buchstab's differential-delay
equation, which concerns large integers $N$ divisible by no prime less than $N^{1/u}$.
Discussion of the latter problem occurred in letters to the journal Nature during the second world
war, from the physicists Frederick Lindemann and Wendell Furry. We recount how Furry
evaluated a dilogarithm in reply to a puzzle resulting from Mertens' third theorem,
raised by Lindemann. We refine Furry's analysis to include multiple 
polylogarithms of  weights up to 200.

\section{Introduction}

In 1930, Karl Dickman considered {\em smooth} numbers, all of whose prime divisors are
smaller than a certain magnitude~\cite{Di}. In 1937, Aleksandr Buchstab considered
{\em rough} numbers, none of whose prime divisors are smaller
than such a magnitude~\cite{Buch}. These cognate problems received
further attention after the second world war, notably from Nicolaas de Bruijn~\cite{dBr1,dBr2,dBr3}.
During that war, Lindemann and Furry addressed the problem of
rough numbers, which clearly include primes.

Replying~\cite{F42} to a letter~\cite{L41a} by Lindemann, 
{\em Number of primes and probability considerations}, Furry deftly computed the asymptotic 
density of numbers of size $N$ divisible by no prime $p\le N^{1/4}$. These comprise
primes, semiprimes and integers with three prime divisors, all of which are greater than $N^{1/4}$.
From this analysis, he obtained an approximation to $e^{-\gamma}$ with an absolute
error less than $1.25\times10^{-6}$, where the Euler--Mascheroni constant $\gamma$ emerges
from a combination of Mertens' third theorem with probabilistic reasoning, as suggested by Lindemann.
In achieving this, Furry dealt with dilogarithms~\cite{Lewin},
in advance of their application to quantum field theory by  Schwinger~\cite{JS},
and with a Mertens product, in advance of Buchstab's later work~\cite{Buch2} in 1951, 
which cites Brun's use of the sieve of Eratosthenes~\cite{Brun}.

The Dickman problem~\cite{Di} for smooth~\cite{PM} numbers
was considered by many post-war authors, both before~\cite{CV,Rama} and after~\cite{BK,Cham,AH,HT,DK}
de Bruijn~\cite{dBr2}, since issues of smoothness and semi-smoothness~\cite{BP} are of importance
in situations where one hopes to make progress in factorizing large composite numbers~\cite{RB}.

In this article, we explain the context of Furry's letter on the Buchstab problem~\cite{Buch}, written 
when he was working
on radar in USA, and Lindemann's letter, written when he was Churchill's scientific advisor,
much preoccupied with the war effort in UK. Moreover
we develop Furry's analysis to include multiple polylogarithms
of weights up to 200, using techniques devised for calculations in quantum field theory,
where a pressing need to understand phenomena at the Large Hadron Collider (LHC)
has resulted in computational progress, albeit for integrals with lesser weights~\cite{NSU,VW} . 

For $u>1$ the Dickman function $\rho(u)$ satisfies the differential-delay equation~\cite{dBr2}
\begin{equation}u\rho^\prime(u)=-\rho(u-1)\label{DE}\end{equation}
with $\rho(u)=1$ for $u\in[0,1]$. To compute it, we define multiple polylogarithms
\begin{equation}M_{j,n}(y)=\sum_{n_1>\ldots>n_j>0}\;\prod_{i=1}^j\frac{z_i^{n_i}}{n_i},\quad
z_1=\frac{y}{n},\quad z_i=\left(\frac{n+2-i}{n+1-i}\right)~\text{for}~i>1\label{Mjn}\end{equation}
with nested sums that converge rapidly for real $y\in[0,1]$ and integers $n>j>0$.

{\bf Theorem~1:} {\em For $y\in[0,1]$ and $n>0$,}
\begin{equation}\rho(n-y)=\rho(n)+\sum_{j=1}^{n-1}\rho(n-j)M_{j,n}(y).\label{Th1}\end{equation}

This will be proved in Section~2, where it leads to a practical method for computing $\rho(u)$,
for real $u\in[1,200]$, thanks to the procedure  {\tt polylogmult} of {\tt Pari/GP}~\cite{Pari},
which implements Henri Cohen's efficient extension~\cite{HC} of an algorithm devised by 
P.\ Akhilesh~\cite{Akh1,Akh2} for Don Zagier's multiple zeta values~\cite{MZV,DZ}. 

For $u>2$, the Buchstab function $\omega(u)$ satisfies the differential-delay equation~\cite{dBr1}
\begin{equation}u\omega^\prime(u)=\omega(u-1)-\omega(u)\label{BE}\end{equation}
with $u\omega(u)=1$ for $u\in[1,2]$. Let $\sigma(u)=(u+1)\omega(u+1)$.

{\bf Theorem~2:}  {\em For real $y\in[0,1]$ and $n>0$,}
\begin{equation}\sigma(n-y)=\sigma(n)+\sum_{j=1}^{n-1}(-1)^j\sigma(n-j)M_{j,n}(y).\label{Th2}\end{equation}

Thanks to this, we shall solve the Buchstab problem in parallel with the Dickman problem,
merely by changing the signs of polylogs with odd weights. For real $u>1$, the results have the form
\begin{equation} \sigma(u)=1+\sum_{u>k>0}P_k(u),\quad\rho(u)=1+\sum_{u>k>0}(-1)^kP_k(u)\label{both}\end{equation}
where the sums terminate at the largest integer less than $u$ and $P_k(u)$ has support only for $u>k$,
where it is positive and of pure polylogarithmic weight $k$.  

For $u\ge k>0$, the boundary condition $P_k(k)=0$
and the differential equation $uP_k^\prime(u)=P_{k-1}(u-1)$ determine $P_k(u)$ as an iterated integral,
on the understanding that $P_0(u)=1$ for $u\ge0$. Thus  $P_1(u)=\log(u)$, for $u\ge1$.
For $u\ge2$, we have
\begin{equation}P_2(u)=\frac{\log^2(u)-\zeta_2}{2}+{\rm Li}_2\left(\frac{1}{u}\right)\label{P2}\end{equation}
where ${\rm Li}_k(z)=\sum_{n>0}z^n/n^k$, for $|z|<1$, is a classical polylogarithm~\cite{Lewin}
with ${\rm Li}_1(z)=-\log(1-z)$ and $z{\rm Li}_k^\prime(z)={\rm Li}_{k-1}(z)$, for $k>1$,
where $\zeta_k={\rm Li}_k(1)$. 
To prove~(\ref{P2}) we use the evaluation~\cite{Lewin}
${\rm Li}_2(\frac12)=\frac12\zeta_2-\frac12\log^2(2)$, which shows that the condition $P_2(2)=0$
is satisfied. The differential equation $uP_2^\prime(u)=\log(u-1)$, for $u\ge2$, is then verified
as the identity $\log(u)+\log(1-1/u)=\log(u-1)$.

In 1941, Lindemann wrote a letter~\cite{L41a} to the editors of the journal Nature, where it was published with
the title {\em Number of primes and probability considerations}. He noted
the  appearance of Euler's constant $\gamma$ in Mertens' third theorem~\cite{M}
\begin{equation}\lim_{x\to\infty}\log(x)\prod_{p\le x}\left(1-\frac{1}{p}\right)=e^{-\gamma},\label{M3}\end{equation}
for products over primes, and its conspicuous absence from the prime number theorem.  
In a reply~\cite{F42} by Furry, we found an intriguing approximation
\begin{equation}4e^{-\gamma}-\log(3)-\log(2)\log\left(\frac32\right)+\frac12\sum_{k=1}^\infty\frac{1}{k^24^k}
\approx1.0000050\label{Fk}\end{equation}
whose probabilistic origin we shall review in Section 3.

In Section~4, we study the behaviour of $P_k(u)$, $\sigma(u)$ and $\rho(u)$ at large $u$.
Section~5 concerns the discrete situation, where we study divisibility properties of integers
in finite ranges. Section~6 offers comments and conclusions.

\section{Algorithms and proofs}

The multiple polylogarithms $P_k(u)$ in the terminating series~(\ref{both}) are $k$-fold iterated integrals that
may be evaluated as nested multiple sums~\cite{BBB3,RV}.

{\bf Theorem~3:} {\em For $u\ge k>0$,}
\begin{equation}
P_k(u)= \sum_{n_1>\ldots>n_k>0}\;\prod_{i=1}^k\frac{z_i^{n_i}}{n_i},\quad
z_i=1-\frac{1}{u-k+i}\label{Pk}.\end{equation}

The proof is given later. Here we remark on the rate of convergence.
The nested sum in~(\ref{Pk}) clearly converges, since $|z_i|<1$.
Yet as $u$ increases, the convergence becomes slower. For $u\in[k,k+1]$,
{\tt polylogmult} is very efficient, giving 50 digits of $P_{10}(11)$ in
less than 3 milliseconds on a 3.1 GHz machine. Yet $P_{10}(13)$ takes  
34 milliseconds. If we ask {\tt Pari/GP} for $P_{10}(14)$ we are politely told:
{\em sorry, polylogmult in this range is not yet implemented.}
To overcome this, we use Theorem~2.

\subsection{Algorithms}

To evaluate $\rho(u)$ and $\sigma(u)$ with $u\in[1,N]$, we begin
by evaluating a triangular matrix of constants, which will be stored.
Consider polynomials in $\lambda$ defined recursively by
\begin{equation}S_n=S_{n-1}-\sum_{j=1}^{n-1}(-\lambda)^jS_{n-j}M_{j,n}(1)\label{rec}\end{equation}
for $n\in[2,N]$, with $S_1=1$. Then Theorem~2 gives $\sigma(n)$ as the value of $S_n$ at $\lambda=1$.

{\bf Algorithm~1:} With symbolic $\lambda$ and numerical $M_{j,n}(1)$, perform recursion~(\ref{rec})
and store the coefficient of $\lambda^k$ in $S_n$ as $P_k(n)$, for $0\le k<n\le N$.

With $y=1$ in the multiple polylogarithms~(\ref{Mjn}) this evaluates, for example,
\begin{align} 
P_1(4)&=M_{1,2}+M_{1,3}+M_{1,4},\label{P14}\\
P_2(4)&=M_{1,2}(M_{1,3}+M_{1,4})+M_{1,3}M_{1,4}-M_{2,3}-M_{2,4},\label{P24}\\
P_3(4)&=M_{1,2}M_{1,3}M_{1,4}-M_{2,3}M_{1,4}-M_{2,4}M_{1,2}+M_{3,4}.\label{P34}
\end{align}

{\bf Algorithm~2:} Let $n$ be the ceiling of $u$. For $k<n$, evaluate
\begin{equation}P_k(u)=\sum_{j=0}^k(-1)^jP_{k-j}(n-j)M_{j,n}(n-u)\label{alg2}\end{equation}
on the understanding that the empty product gives $M_{0,n}(y)=1$ at $j=0$.

\subsection{Small weights}

There is a limitation to these algorithms at large $u$. The poor convergence in Theorem~3,
for increasing $u-k$, has been traded for an increasing  number of products in Algorithm~2,
with constants from Algorithm~1 that involve products of products. This situation is well
understood in the application of perturbative quantum field theory to experimental
high-energy physics. The formal parameter $\lambda$ of Algorithm~1 is eventually set
to $\lambda=1$, to solve the Buchstab equation, or to $\lambda=-1$, to solve the
Dickman equation. It mimics the coupling constant of a perturbative expansion.
Higher powers of $\lambda$ involve an increase of the maximum weight by 1.
In high-energy physics, the weight may increase by 2
at the next order in $\lambda$, if the result is reducible to multiple polylogarithms~\cite{BS,SYM}.

The Dickman argument $u$ is akin to a kinematic variable in high-energy physics,
where there are situations in which polylogs become large in certain regions
of phase space. In such a case, it is useful to have a leading-logarithm approximation~\cite{SYM}.
In the Dickman problem, the leading power of a log is easily determined: $P_k(u)$ is asymptotic
to $\log^k(u)/k!$. On the basis of high-precision results for $k<10$, it was conjectured in~\cite{DB}
that at large $u\gg k$
\begin{equation}
P_k(u)=\sum_{j=0}^k\frac{D_{k-j}\log^j(u)}{j!}+o(1),\quad 
\sum_{k=0}^\infty D_kz^k=D(z)=\frac{e^{-\gamma z}}{\Gamma(1+z)}\label{db}\end{equation}
and this was later proved in~\cite{Sound}. 

This behaviour is not immediately apparent from the algorithms, which make
calls to {\tt polylogmult} via auxiliaries in~(\ref{Mjn}) with good convergence.
However, we must combine many products of such terms when $u$ is large. 
Consider weight 2, where we know that $P_2(u)=\frac12\log^2(u)-\frac12\zeta_2+{\rm Li}_2(1/u)$, for $u\ge2$. 
This form is optimal at large $u$, requiring a single call to {\tt polylog} with a small argument $1/u$. 
In contrast, our method appears profligate at large $u$. Let $N$ be the ceiling of $u$.
Then {\tt polylogmult} evaluates $P_2(n)=P_2(n-1)+M_{2,n}(1)-\log(n-1){\rm Li}_1(1/n)$, for $n\in[3,N]$.
Finally, we obtain $P_2(u)=P_2(N)-M_{2,N}(N-u)+\log(N-1)\log(N/u)$.
This involves roughly  $10N$ times more work than is needed for a single value of $P_2(u)$, with $u$ close to $N$,
because we have taken about $N$ steps and incurred a slowdown by a factor of about  10
by using  {\tt polylogmult} instead of the highly tuned classical {\tt polylog} routine.
Yet the algorithms have merits: they work at all weights and we store the constants from Algorithm~1,
which are reused in various applications of Algorithm~2. 

Now consider weight 3, where efficient reduction to classical polylogs is also possible.

{\bf Theorem~4:} {\em For $u\ge3$,}
\begin{align}P_3(u)&=\tfrac12{\rm Li_3}\left(\frac{1}{u(2-u)}\right)-{\rm Li}_3\left(\frac1u\right)
-{\rm Li}_3\left(\frac{1}{2-u}\right)+{\rm Li}_2\left(\frac1u\right)\log(u-2)\nonumber\\{}&
+\tfrac13\zeta_3-\tfrac12\zeta_2\log(u)
+\tfrac{1}{12}\log^3(u(u-2))-\tfrac12\log^2(u-2)\log(u).\label{Th4}\end{align}

As usual, we postpone the proof and concentrate on the structure, which is optimal. 
The 4 classical polylogs in~(\ref{Th4}) converge very rapidly at large $u$. 
The remaining terms give $P_3(u)=\frac13\zeta_3-\frac12\zeta_2\log(u)+\frac16\log^3(u)+o(1)$ at large $u$,
in accord with~(\ref{db}). For large $u$ of size $N$, high precision
evaluation of~(\ref{Th4}) using {\tt polylog} is about $2N$ times faster than an $N$-step 
process using {\tt polylogmult}.

We expect to need {\tt polylogmult} at weight 4. Yet there is still an efficient way
to evaluate $P_4(u)$, without the $N$-step process. For $\frac12\ge y>0$, we define
\begin{align}
E_4(y)&=\tfrac{19}{16}\zeta_4
-3\,{\rm Li}_4(-y)
+3\,{\rm Li}_3(-y)\log(y)
-\tfrac32\,{\rm Li}_2(-y)\log^2(y)\nonumber\\{}&
+\tfrac12\,{\rm Li}_1(-y)\log^3(y)
+\tfrac18\log^4(y)
\label{E4}\\
H_4(y)&=\int_0^y\left(\log\left(\frac{x}{1+2x}\right)\,{\rm Li}_2(x)
+\frac12\log^2(x)\,{\rm Li}_1(-2x)\right)\frac{{\rm d}x}{x(1+x)}\label{H4}
\end{align}
with $E_4(y)$ easy to compute, 
while $H_4(y)$ is harder, yet vanishes as $y\to0$, since $H_4(y)=O(y\log^2(y))$.

{\bf Theorem~5:} {\em With $u\ge4$ and $y=1/(u-2)$,}
\begin{equation}P_4(u)=P_3(u-1)\log(u)
+\frac{\zeta_2}{4}\left(2\,{\rm Li}_2\left(1-\frac1u\right)+\log^2(u)\right)
-E_4(y)-H_4(y).\label{Th5}\end{equation}

The proof is given later. The important point  here is that the hard integral $H_4(y)$
becomes easier at large $u$ and hence small $y=1/(u-2)$, while for $u$ close to 4 we simply use Theorem~3.
There is a critical value $u_c\in[6.0,\,7.0]$ below which Theorem~3 should be used for $P_4(u)$ and above which 
Theorem~5 should be used. We are grateful to Steven Charlton for automating the small $y$ expansion 
of $H_4(y)$ by reduction to multiple polylogarithms of the general {\tt polylogmult} form
\begin{equation}{\rm Li}_{s_1,...,s_d}(z_1,...,z_d)
=\sum_{n_1>\ldots>n_d>0}\;\prod_{i=1}^d\frac{z_i^{n_i}}{n_i^{s_i}}\quad
\label{Lgen}\end{equation}
with depth $d$ and weight $k=\sum_{i=1}^d s_i$. The result is
\begin{align}
H_4(y)&=\tfrac12C_2(y)\log^2(y)+C_3(y)\log(y)+C_4(y)\label{H4e}\\
C_2(y)&={\rm Li}_2(-2y)+{\rm Li}_{1,1}(-y,2)\label{C2}\\
C_3(y)&={\rm Li}_{3}(y)
-{\rm Li}_{3}(-2y)
+{\rm Li}_{1,2}(-y,-1)
-{\rm Li}_{2,1}(-y,2)\label{C3}\\
C_4(y)&={\rm Li}_{4}(-2y)
-{\rm Li}_{4}(y)
+{\rm Li}_{2,2}(-2y,-\tfrac12)
-{\rm Li}_{2,2}(-y,-1)\nonumber\\&{}
+{\rm Li}_{3,1}(-2y,-\tfrac12)
+{\rm Li}_{3,1}(-y,2)
+{\rm Li}_{3,1}(y,2)
+{\rm Li}_{1,1,2}(-y,2,-\tfrac12)\nonumber\\&{}
+{\rm Li}_{1,2,1}(-y,-1,-2)
+{\rm Li}_{1,2,1}(-y,2,-\tfrac12)
\label{C4}\end{align}
with ordering of subscripts and arguments as in~\cite{BBB3}, which is reversed in~\cite{SCh}.
This evaluates the Furry probability $P_4(u)$ at 50-digit precision
in about 4 milliseconds, for any real value $u>4$. It takes less than $0.6$ seconds
to achieve 1000-digit precision. 

This method may be continued to weights $k>4$, by defining a set of
auxiliary functions, $F_k(u)$, with $F_0(u)=1$, for $u\ge0$, and $F_k(k)=0$ for $k>0$.
These solve $(u-k)F_{k+1}^\prime(u)=F_k(u-1)$ for $u+1>k\ge0$.
It follows that $F_1(u)=P_1(u)=\log(u)$, for $u>1$.
For $u>k>1$, we use the recursion~\cite{DB}
\begin{equation}F_k(u) =\sum_{j=0}^{k-1}(-1)^{k-j-1}P_{k-j}(u-j)F_j(u).\label{frec}\end{equation}

Inspection of~(\ref{alg2}) shows that Algorithm~2 relies on an algebraic structure
similar to that in~(\ref{frec}).
Algorithm~1 packs products of auxiliary constants,
$M_{j,n}(1)$, to furnish constants of interest, namely $P_k(n)$ at integer arguments $n>k$,
which we prudently store, at high precision.
Algorithm~2 efficiently uses those stored constants, to evaluate
further values of actual interest, namely the Furry probabilities $P_k(u)$ for real $u>k$, 
by evaluating values of $M_{j,n}(y)$, with $y=n-u\in(0,1)$ and $j\in[1,k]$.
The following theorem~shows the utility of the further packing of products in~(\ref{frec}).

{\bf Theorem~6:} {\em For $u\ge k>n\ge0$,} 
\begin{equation}I_{k,n}(u)=\int_k^u P_{k-n-1}(x-n-1)\frac{F_n(x)\,{\rm d}x}{x-n}
=\sum_{j=0}^{n}(-1)^{n-j}P_{k-j}(u-j)F_j(u).\label{Th6}\end{equation}

At $n=0$, this yields an obvious result, namely that $I_{k,0}(u)=P_k(u)$ is the integral 
that solves $uP^\prime_k(u)=P_{k-1}(u-1)$ with $P_k(k)=0$. Less trivially, at $n=k-1$,
Theorem~6 proves the recursion~(\ref{frec}). Most usefully, it shows that {\em any} evaluation
of an integral $I_{k,n}(u)$,  with $k>n\ge0$,  
will serve our avowed purpose of evaluating the Furry probability $P_k(u)$. 
It was by this means that we derived~(\ref{Th5}), in terms of integrals over a product 
of a dilogarithm and a logarithm, or products of 3 logs.

More generally, at weight  $k<10$, it suffices to integrate over a sum of products of polyogs whose 
individual weights do not exceed $k/2$. Moreover, the asymptotic behaviour of 
the integrand is determined by~(\ref{db}) and~(\ref{frec}), enabling one to separate the
integral into an easier and a harder part, as in~(\ref{E4},\ref{H4}), with classical polylogs in
the easier part, which dominates at large $u$, while the harder part vanishes as $u\to\infty$
and is not needed for $u-k<3$, where Theorem~3 may be used.
By this means it was possible to obtain high-precision values~\cite{DB} of the Dickman constants
$D_k$, for weights $k<10$, and to infer their generating function $D(z)$ in~(\ref{db}).
For example, the exact value of
\begin{equation}D_9=\tfrac19\zeta_9-\tfrac{1}{14}\zeta_7\zeta_2+\tfrac{1}{80}\zeta_5\zeta_4
-\tfrac{5}{384}\zeta_3\zeta_6+\tfrac{1}{162}\zeta_3^3\,\approx\,0.0016850\label{D9}\end{equation}
was discovered  after numerical  quadrature for an integral, over $x\in[9,\infty]$, whose integrand
was obtained by subtracting a polynomial of  degree 8 in $\log(x)$ from the product
$P_4(x-5)F_4(x)$.
An $N$-step method would be of no help here. With $k=9$ and $u=200$, the $o(1)$ term in~(\ref{db})
is approximately $0.039906$. While small, compared with $P_9(200)\approx1.6226$,
it is an order of magnitude greater than $D_9$.

\subsection{Efficiency and complexity}

The algorithms follow good banking practice. We incur a debt
in Algorithm~1 that is amortized by making many calls to the faster Algorithm~2.

By way of example, consider $P_k(u)$ with weights $k<10$ and real arguments $u<201$. Discarding terms
of order $\lambda^{10}$ in~(\ref{rec}), the 200 steps of Algorithm~1 require 1764 evaluations of~(\ref{Mjn}) at $y=1$,
taking about $3.3$ seconds, at 100-digit working precision, and about 190 seconds at 1000-digit precision, when
calling {\tt polylogmult} in a single thread on a modest $3.1$~GHz processor. Then Algorithm 2 requires
only 9 calls with $y\in(0,1)$, to evaluate $P_k(u)$ for all $k<10$
and fixed real $u=n-y<201$, where $n=\lceil u\rceil$ is the ceiling of $u$. This takes about $1.1$ seconds
at 1000-digit precision for a random real value of $u<201$. The time taken by Algorithm~2 is commensurate with the 
methods of Theorems~4 and~5, which evaluate 1000 digits of $P_k(u)$ with $k<5$ in less than $0.7$ seconds.
For  $k\in[5,9]$, Algorithm~2 is more efficient than numerical quadrature based on Theorem~6,
if one discounts the debt accrued by Algorithm~1.

For $N\ge n>j\gg1$, the complexity of a single call to {\tt polylogmult} for $M_{j,n}(1)$, in 
the $N$-step Algorithm~1, increases faster than $j^2b^2$, to achieve an absolute error 
less than $1/2^b$, since it requires  $O(j^2b)$ floating point operations on numbers with $O(b)$ bits. 
Thus the complexity of Algorithm~1 increases faster than $N^4D^2$, at a working precision of
$D$ decimal digits. An initial investigation with $N=101$ took less than an hour, in a single thread
with a working precision of $D=350$ decimal digits, which was ample for the work in Section~4  
at weights $k\le100$.  Thereafter we stepped up to $N=201$ at $D=1000$ digits. This took less than
a CPU-week. The effort is embarrassingly parallelizable, since it consists in the 
accumulation of about $\frac12N^2$ independent constants, each at a cost of order $N^2D^{2+\epsilon}$,
with $\epsilon=o(1)$ for Sch\"onhage–Strassen multiplication, or $\epsilon\approx0.465$
for Toom-Cook multiplication~\cite{CP}.

\subsection{Proofs by descent in weight}

Proofs of the first three theorems rely on the relation between 
nested sums and iterated integrals developed in~\cite{BBB1,BBB2,BBB3,ABG,MUW}.
Consider the $k$-fold iterated integral 
\begin{equation}G_k({\bf a},y)=\int_0^y\frac{{\rm d}x_1}{x_1-a_1}
\int_0^{x_1}\frac{{\rm d}x_2}{x_2-a_2}\ldots\int_0^{x_{k-1}}\frac{{\rm d}x_k}{x_k-a_k}\label{Gk}\end{equation}
where {\bf a} is a $k$-letter word denoting concatenation of the constants $a_i$, all of which are greater than $y$, for the
cases needed here. 
Differentiating with respect to $y$, we obtain  $(y-a_1)G_k({\bf a},y)=G_{k-1}({\bf b},y)$
where {\bf b} is obtained by removing the first letter of {\bf a}, on the understanding that $G_0(\text{\o},y)=1$
where {\o} is the empty word. Expanding integrands in $x_i/a_i$, we obtain
\begin{equation}(-1)^kG_k({\bf a},y)=\sum_{n_1>\ldots>n_k>0}\;\prod_{i=1}^k\frac{z_i^{n_i}}{n_i},\quad
z_1=\frac{y}{a_1},\quad z_i=\frac{a_{i-1}}{a_i}~\text{for}~i>1.\label{dict}\end{equation}

{\em Proof of Theorem~1:} Applying the dictionary in~(\ref{dict}) to~(\ref{Mjn}), 
we obtain $M_{j,n}(y)=(-1)^jG_j({\bf a},y)$ with $a_i=n+1-i$ and $n>j\ge i>0$. 
Hence $(n-y)M_{j,n}^\prime(y)=M_{j-1,n-1}(y)$, on the understanding
that the empty product is unity. Observing that~(\ref{Th1})
is true at $y=0$, for all $n>0$, and also true at $n=1$, for all $y\in[0,1]$,
we assume that $n>1$, differentiate and multiply by $(n-y)\ge1$, to obtain
\begin{equation}\rho(n-1-y)=\rho(n-1)+\sum_{j=2}^{n-1}\rho(n-j)M_{j-1,n-1}(y)\label{Pr1}\end{equation}
with the Dickman equation~(\ref{DE}) used on the left and the empty product separated on the right.
Observing that~(\ref{Pr1})
is equivalent to~(\ref{Th1}), with $n$ replaced by $n-1$, we prove the latter for all $n>0$,
by the method of descent. If~(\ref{Th1}) were false for some integer $n>1$, then it would be 
false for every positive integer less than $n$, which is impossible, since~(\ref{Th1}) is true for $n=1$.~$\square$

{\em Proof of Theorem~2:} The Buchstab problem requires that $\sigma(u)=(u+1)\omega(u+1)$
solves $u\sigma^\prime(u)=\sigma(u-1)$ for $u>1$, with $\sigma(u)=1$ for $u\in[0,1]$. We proceed as in Theorem~1,
obtaining
\begin{equation}\sigma(n-1-y)=\sigma(n-1)-\sum_{j=2}^{n-1}(-1)^j\sigma(n-j)M_{j-1,n-1}(y)\label{Pr2}\end{equation}
from differentiation of~(\ref{Th2}) with $n>1$. Observing that~(\ref{Pr2})
is equivalent to~(\ref{Th2}), with $n$ replaced by $n-1$, we prove the latter for all $n>0$.~$\square$

{\em Proof of Theorem~3:} Let $R_k(y)=P_k(k+y)$. Then $(k+y)R_k^\prime(y)=R_{k-1}(y)$
for $y>0$ and $k>0$, on the understanding that $R_0(y)=1$. Hence $R_k(y)=G_k(\tilde{\bf a},y)$
with $\tilde{a}_i=-(k+1-i)$. Now make the transformation $x_i\to y-x_i$, to obtain 
$R_k(y)=(-1)^kG_k({\bf a},y)$ where {\bf a} is obtained from $\tilde{\bf a}$ by reversal and subtraction, 
giving $a_i=y-\tilde{a}_{k+1-i}=y+i$. Then the dictionary in~(\ref{dict}) gives~(\ref{Pk}), with
$z_i=(y+i-1)/(y+i)=1-1/(u-k+i)$.~$\square$

To prove Theorems 4 and 5, we also use descent in weight.
The general method is as follows.
Suppose that we wish to prove a linear relation between terms of the same weight,
as in~(\ref{Th4}), where each term has a rational coefficient multiplying a polylog, or product of polylogs,
whose arguments are rational functions of a variable $u$.
\begin{enumerate}
\item Prove that there is a value of $u$ for which the claim holds.
\item Differentiate with respect to $u$, to obtain a relation between terms of lesser weight, 
with factors that depend rationally on $u$.
\item Take partial fractions to arrive at a set of claims of the same type as before. 
\item For each of these proceed as before until arriving at weight 1, where rational 
linear relations between logs are easy to prove.
\end{enumerate}
Leonard Lewin, following in the footsteps of 
William Spence (1777--1815) and Ernst Kummer (1810--1893), 
was an able practitioner of this art~\cite{Lewin}.

{\em Proof of Theorem~4:} Claim~(\ref{Th4}) agrees with~(\ref{db}) as $u\to\infty$.
The differential equation $P^\prime_3(u+1)=P_2(u)/(u+1)$
requires that $A_2(u)/(u-1)+B_2(u)/(u+1)=0$, for $u>2$, with weight 2 numerators
of partial fractions given by
\begin{align}
A_2(u)&=
{\rm Li}_2\left(\frac{1}{u+1}\right)
+{\rm Li}_2\left(\frac{1}{1-u}\right)
-\tfrac12{\rm Li}_2\left(\frac{1}{1-u^2}\right)
+\tfrac14{\rm Li}_1^2\left(\frac{2}{u+1}\right)\label{T1}\\
B_2(u)&=A_2(u)-{\rm Li}_2\left(\frac{1}{u}\right)-{\rm Li}_2\left(\frac{1}{1-u}\right)
-\tfrac12{\rm Li}_1^2\left(\frac{1}{u}\right)\label{T2}\end{align}
each of which vanishes as $u\to\infty$ and has  a vanishing differential at weight 1 .~$\square$

We can now prove the notable identity
\begin{equation}2\,{\rm Li}_3\left(\frac13\right)-{\rm Li}_3\left(-\frac13\right)=\frac{13\,\zeta_3-\pi^2\log(3)+\log^3(3)}{6}
\label{id3}\end{equation}
by setting $u=3$ in~(\ref{Th4}).
Such relations between constants, rather than functions~\cite{AD}, are often hard to prove. Identity~(\ref{id3})
may be obtained, with ingenuity,  from a special case of a bivariate result  found by Spence, rediscovered by Kummer and
cleaned up by Lewin. This reduces a linear combination of 10 trilogarithms to products. 
With $x=-1$ and $y=\frac13$ in Equation 6.107 of~\cite{Lewin}, it degenerates to
a reduction of 5 trilogarithmic constants  to products,
\begin{equation}
{\rm Li}_3(-3)
-6\,{\rm Li}_3(-1)
-6\,{\rm Li}_3\left(\frac13\right)
+2\,{\rm Li}_3\left(-\frac13\right)
+2\zeta_3=
\frac{\pi^2\log(3)-2\log^3(3)}{3}.\label{SK}\end{equation}
This agrees with~(\ref{id3}) after eliminating
${\rm Li}_3(-3)={\rm Li}_3(-\frac13)-\frac16\pi^2\log(3)-\frac16\log^3(3)$ and
${\rm Li}_3(-1)=-\frac34\zeta_3$. Theorem~4 avoids such mental gymnastics.

{\em Proof of Theorem~5:} The constant $\frac{19}{16}\zeta_4$ in~(\ref{E4}) makes~(\ref{Th5}) consistent 
with~(\ref{db}) as $u=2+1/y\to\infty$. We use the differential equations
$P_4^\prime(u)=P_3(u-1)/u$ and $P_3^\prime(u-1)=P_2(u-2)/(u-1)$ to obtain 
\begin{equation}2yE^\prime_4(y)=
\frac{2\,{\rm Li}_2(y)+\log^2(y)}{1+y}\log\left(\frac{y}{1+2y}\right)
-2yH_4^\prime(y)=\frac{\log^3(y)}{1+y}\label{w3}
\end{equation}
at weight 3. This is verified by straightforward differentiation of $E_4(y)$ in~(\ref{E4}).~$\square$

Consequently, we prove the integer relation
\begin{gather}16H_4(\tfrac12)=
\zeta_4+48\,{\rm Li}_4(-\tfrac12)-8\zeta_2{\rm Li}_2(\tfrac14)
+(12\,{\rm Li}_3(\tfrac14)+16\zeta_2\log(3)-42\zeta_3)\log(2)\nonumber\\{}
+(12\,{\rm Li}_2(\tfrac14)-4\zeta_2)\log^2(2)
-8\log(3)\log^3(2)+10\log^4(2)\label{H4h}\end{gather}
and evaluate the hard integral~(\ref{H4}) at the upper limit $y=\frac12$, where
$P_4(4)=0$. As remarked of~(\ref{id3}), such integer
relations between constants are harder to prove than to
discover empirically~\cite{PSLQ,BB} using {\tt PSLQ}, or {\tt lindep} in {\tt Pari/GP}.

{\em Proof of Theorem~6:} At $n=0$, we obtain $I_{k,0}(u)=P_k(u)$, since $F_0(x)=1$ and $P^\prime_k(x)=P_{k-1}(x-1)/x$.
Both sides vanish at $u=k$, so there is no constant of integration.
Since $F_n(x)/(x-n)$ is the differential of $F_{n+1}(x)$, we may integrate by parts, for $k>n+1>1$, obtaining 
$I_{k,n}(u)=P_{k-n-1}(u-n-1)F_{n+1}(u)-I_{k,n+1}(u)$, after using 
$P^\prime_{k-n-1}(x-n-1)=P_{k-n-2}(x-n-2)/(x-n-1)$.
Again there is no constant of integration. 
Then~(\ref{Th6}) follows by induction on $n\in[0,k-1]$, with $I_{k,k-1}(u)=F_k(u)$ proving
recursion~(\ref{frec}).~$\square$

\section{Six authors in search of a chronicle}

This section is for readers who have an interest in
mathematics as a human activity, rather than an abstract body of knowledge.
It concerns 6 authors linked by the subject of the Lindemann--Furry letters 
to the journal Nature. The names of two of these 6 remain unknown,
despite our earnest efforts to find out who they were.

\subsection{Aleksandr Adolfovich Buchstab (1905--1990)}

Buchstab gained his doctorate in 1939
from Moscow State University, advised by Kinchin.
The portal {\tt Math-Net.Ru}
lists 8 of his publications~\cite{MNR} in the period 1933--1967.
All of these are
single-author articles, in Russian, dealing with number theory.
His 1937 paper~\cite{Buch}
has an abstract in German and gives his affiliation at the time as Baku State University
in Azerbaijan.

Buchstab gives the expansion of $\sigma(u-1)$ 
for $u>2$ as a terminating series of iterated integrals, equivalent to that in~(\ref{both}). 
The definition $\omega(u)=\sigma(u-1)/u$ came 13 years later, from de Bruijn~\cite{dBr1}, 
who references~\cite{Buch}, adding a footnote saying that his own 
article can be read independently from Buchstab's. 
It is unclear whether Furry might have known about Buchstab's article when writing 
from Harvard to Nature in 1942.
It seems quite likely that Lindemann did not know of it.\

\subsection{Frederick Alexander Lindemann (1886--1957)} 

Lindemann gained his doctorate in 1911 from the Friedrich Wilhelm University in Berlin,
advised by Nernst. Returning to UK, he served at the Royal Aircraft Factory during the first world war.
He became director of the Clarendon Laboratory in Oxford in 1919.
As noted by Wright~\cite{EW}, Lindemann had
an active interest in the theory of numbers. His proof~\cite{L33}
of the fundamental theorem of arithmetic,
by Fermat's method of descent, was referenced by Hardy and Wright~\cite{HW}
for its simplicity and elegance. 

Lindemann's friendship with Winston Churchill
resulted in an appointment as chief scientific advisor, when Churchill succeeded Chamberlain 
as prime minister in 1940. From 6 June 1941, Lindemann was often referred to by
his new title, Lord Cherwell. Somewhat symmetrically, the title Fellow of the Royal Society
(FRS) was conferred on Churchill in the cabinet room on 12 June 1941.

Notwithstanding intense involvement with the war effort,
Lindemann had the habit, at weekends, of travelling in his chauffeur-driven Rolls Royce
from London to Oxford, where he could relax in his college rooms, with a bottle of champagne, a popular
illustrated magazine and the Quarterly Journal of Mathematics~\cite{Farm}.
One is tempted to imagine
that an understandable wartime concern with applied statistics and probabilities, combined with
his interest in number theory, led to the letter written in Oxford on Saturday 13 September 1941 
and published~\cite{L41a} in Nature, on 11 October, with the title {\em Number of primes and probability considerations}.

His argument, in paraphrase, is as follows. The probability that a large random number $N$ is not divisible
by a prime $p$ is $(1-1/p)$. Now consider all the primes $p\le x$ and suppose that
these probabilities are independent. Then Mertens'  third theorem gives $e^{-\gamma}/\log(x)$
as the product of probabilities, for large $x$. Now set $x=\sqrt{N}$, to give the sieve of Eratosthenes.
By the prime number theorem, the result should be $1/\log(N)$. Yet we obtain $C/\log(N)$,
with $C=2e^{-\gamma}\approx1.1229$. It follows that the probabilities are not independent.
He adds the comment that ``there is a slight tendency for factors to avoid each other''
and hopes ``that a reader of Nature can throw light on this issue.''

\subsection{John Burdon Sanderson Haldane (1892--1964)}

Haldane had a reputation~\cite{Hnotes} 
for being both brilliant and arrogant. He gained two first-class degrees at Oxford,
one in mathematics, the other in classics and philosophy. His work on 
genetics, evolution, statistics, biochemistry, physiology, ethology, cosmology
and the origins of life was so wide ranging that it took 6 scholars~\cite{Hnotes}
to review it for the Royal Society, on his death.

Yet his reply~\cite{H41} to Nature
was irrelevant to the mathematical question in hand. Haldane argues that
because there is an infinity of positive integers, ``there is no way of choosing one at random'' 
and ends by saying that ``we
biometrists have our difficulties, but at least the number
of men, or even of bacteria, is finite, so biometric
sampling theory can be given a comparatively secure logical
basis." 

To this, Lindemann correctly replied, in a second letter~\cite{L41b},
that ``the discrepancy to which I directed attention can be derived
perfectly well by choosing a number at random from a finite
class."

\subsection{An anonymous Free French Scientist}

The next reply~\cite{FFS41} came from a ``Free French Scientist'' (FFS),
despite the journal's firm statement that no anonymous submission would
be accepted.

Like Lindemann, FFS takes a Mertens product of $(1 - 1/p)$, from $p = 2$ up to $p = p_n$, 
the $n$-th prime, and claims that this is ``valuable only for the whole
interval" of numbers $[p_n+1,\,p_{n+1}^2]$. 
In support of this claim, FFS provides a table up to $n=24$ and claims agreement
to within 1\% between the estimated and actual number of primes for $n\in[17,24]$.
There are errors in the table, but the claim of 1\% agreement for $n=24$ is correct.
There are 1139 primes in the interval $[90,97^2]$, while the FFS prediction rounds to 1133.
(FFS says there are 1129 primes and estimates 1136.)

FFS disputes Lindemann's claim of a discrepancy of more than 10\%, for very large numbers,
and asserts that ``the probabilities are strictly independent and it cannot be any question of a
mysterious tendency of the factors to avoid each other."

Lindemann assumes that FFS is male and replies~\cite{L41b}, saying that
``the agreement he finds is unfortunately spurious\ldots. If he
had proceeded to larger numbers the agreement would vanish
and he would find a discrepancy just as serious as that to
which I directed attention."

A discrepancy is indeed more apparent for larger numbers. For example,
it is more than 3\% for $n=167$, with 77893 primes in the interval $[992, 997^2]$,
while the FFS estimate rounds to 80481.

\subsection{Wendell Hinkle Furry (1907--1984) }

Furry gained  his doctorate from the
University of Illinois at Urbana-Champaign in 1932, with a thesis on
the lithium molecule. His name is familiar to workers in quantum electrodynamics,
thanks to a theorem~\cite{F33} which tells them that they do not need to
compute Feynman diagrams with electron loops coupling to an odd number of photons,
since these cannot contribute to a physical process. (This is not the case for
quark loops coupled to gluons in the non-abelian theory of quantum chromodynamics.)
In 1942, Furry was working on radar at MIT's Rad Lab. 

Furry's letter~\cite{F42} is succinct and deft in its computation. Here we restate his argument, in our own notation.

Consider a range of large numbers of size $N$. Select  those that have no prime divisor $p$
with $p^{u+1}<N$ and $u>1$. Some will be prime. Others will be the product of $k+1$ primes
with $u>k>0$. Let $P_k(u)$ be the asymptotic ratio of those with $k+1$ prime divisors to those that are prime.
Then $P_1(u)=\log(u)$. For $u>2$ we may compute $P_2(u)$ as a dilogarithm. For $4\ge u>3$, 
we need to integrate a dilogarithm to find $P_3(u)$. As $u$ increases, the Mertens approximation
\begin{equation}(u+1)e^{-\gamma}-\sum_{u>k>0}P_k(u)\approx1\label{appr}\end{equation}
becomes dramatically better. Lindemann observed that $2e^{-\gamma}\approx1.1229$.
Furry observed that $3e^{-\gamma}-\log(2)\approx0.99123$. With $u=3$ he obtained
\begin{equation}4e^{-\gamma}-\log(3)-\log(2)\log\left(\frac32\right)+\frac12\sum_{k=1}^\infty\frac{1}{k^24^k}
\approx1.0000050\label{Fk1}\end{equation}
which we give in the form that appears in~\cite{F42}. It was this notable approximation that led us to 
connect the multiple polylogarithms for smooth numbers in the Dickman problem~\cite{DB} to those
for rough numbers in the Lindemann problem~\cite{L41a}.

In~\cite{F42}, Furry was rather smart. At $u=3$, we encounter ${\rm Li}_2(\frac13)$ in~(\ref{P2}),
while~(\ref{Fk1}) is free of $\pi^2$ and involves the even faster converging series 
${\rm Li}_2(\frac14)$. To obtain the latter, one may use 3 tricks that subsequent practitioners 
have learnt from Lewin's valuable book~\cite{Lewin}. 
First use the transformation $z\to-z/(1-z)$, to get from $z=\frac13$ to $z=-\frac12$.
Then use ${\rm Li}_2(-\frac12)+{\rm Li}_2(\frac12)=\frac12{\rm Li}_2(\frac14)$, to get to $z=\frac14$. Finally,
reduce ${\rm Li}_2(\frac12)$ to $\pi^2$ and $\log^2(2)$, as is done by the requirement $P_2(2)=0$. 

Lindemann's reply~\cite{L42} shows that he too was impressed:
``I am grateful to Prof.\ Furry for his sympathy with the outlook of 
a mere physicist and for the light he has thrown on this matter.
Prof.\ E.\ M.\ Wright, some months ago, sent me privately a proof on
somewhat similar lines that the probabilities could not be independent\ldots.
Prof.\ Furry has carried the process a step further and shown the degree
of divergence.''

\subsection{An unidentified colleague of a master bomb maker}

The sixth author in this chronicle is as mysterious as FFS and will be referred to as BM,
since she or he may have been associated with bomb making, supervised by
Millis Rowland Jefferis (1899--1963). We learnt of a manuscript  by BM via entry G.306
in the catalogue of Lindemann's papers held by Nuffield College, Oxford:
``Includes typescript on `Density of Prime Numbers' sent to Cherwell by Jefferis and typescript 
with diagrams and photographs on `Use of Soft Nosed Bombs for Attack of Capital Ships', no author, no date.''
The librarian, Clare Kavanagh, kindly supplied us with a copy of BM's document, which concludes
with a manuscript equation in which an integral appears. We have simplified this, reducing BM's equation to
\begin{equation}K - 1 = \int_{\frac13}^{\frac12}\frac{{\rm d}t}{t}\left(\frac{1}{1-t}-\frac{K}{2t}\right),\label{BM}\end{equation}
where $K$ is not otherwise defined in BM's document.
It is easy so show that~(\ref{BM}) gives $3K/2=1+\log(2)$ and hence $K\approx1.1288$. Yet BM
approximates an integral by a truncated series and writes $1.126$, as a rough value for $K$. 
BM seems pleased that $K$ is close to Lindemann's $C=2e^{-\gamma}\approx1.1229$. We conclude that
BM was undertaking the first step of Furry's development, which gives $1+\log(2)\approx1.6931$ as a fair approximation to
$3e^{-\gamma}\approx1.6844$.

\subsection{Who were FFS and BM?}

BM's understanding of the Mertens discrepancy was better than that shown by FFS,
while  Furry was more adroit in analysis than BM. Accordingly, we assume that these authors are 3 distinct people.

Since  FFS wrote to Nature within two days of the publication~\cite{L41a} in London of Lindemann's first letter, 
we assume that FFS was working in England, in conjunction with Free French forces organized in 
exile by Charles de Gaulle (1890--1970). 
This is consistent with an obituary notice~\cite{FFS42} in Nature,
published on 7 February 1942, where FFS reports that ``the French physicist Fernand 
Holweck died recently in Paris, but the circumstances 
of his death are not clear.'' It was later established that Holweck had
been working with the French resistance, in conjunction with the  Special Operations Executive, SOE,
and had been betrayed, arrested, tortured and killed. For a vivid
account of SOE and its dependence on encryption, see~\cite{SC}. 

Noting that BM's document was sent by Jefferis to Lindemann in an undated
set of papers including both {\em Density of prime numbers} and notes
on bomb making, we assume that BM was associated, in some way,
with MD1, popularly known as ``Winston Churchill’s Toyshop''~\cite{MD1}.
This department designed bombs and also provided sabotage devices for SOE. 
Churchill and Lindemann provided special funds for MD1 and ensured that it
functioned independently of the Ministry of Supply, upon whose remit it encroached.
Jefferis was a senior officer closely involved with MD1, as was Robert Stuart Macrae,
who later wrote a book about its work~\cite{RSM}. 

It is probable that BM and FFS were among a group of scientists and mathematicians 
operating in England within a set of departments that included MD1, SOE, Secret Intelligence (MI6),
where Edward Maitland Wright (1906--2005) worked~\cite{EMW}, 
and the Government Code and Cypher School (GC\&CS) at Bletchley Park,
where Alan Mathison Turing (1912--1954) worked~\cite{AT}. These were of special interest to 
Churchill and his chief scientific advisor, Lindemann, whose reputation
was evidently sufficient for the editors of Nature, Lionel Brimble (1904--1965)
and Arthur Gale (1895--1978), to set aside their rule
against anonymous publication, in the case of FFS, and for
BM to seek private communication, via Jefferis, on the Mertens puzzle.

The style of BM's document follows that of a rough first draft for a
mathematical article, with vocabulary and notation typical of what might be expected
from an academic used to writing in English. Perhaps BM was Wright, though we would like to 
believe that Hardy's co-author was capable of solving~(\ref{BM}) analytically, rather
than by tortuous approximation. In any case, there were many able mathematicians,
including Dennis Babbage (1909--1991), David Champernowne (1912--2000), 
William Tutte (1917--2002), Henry Whitehead  (1904--1960) and Shaun Wylie (1913--2009),
engaged, along with Turing and Wright, on secret war work. In the absence of further
clues, we leave the identity of BM undecided.

It seems clear that FFS was a French speaker, using ``valuable'' to mean
valid (in French, {\em valable}). There were many exiled Free French scientists and intellectuals
working in England at the time. Reading about the activities  of 
Joseph Cathala (1892--1969), 
Jules Gu\'eron (1907--1990),
Hans von Halban (1908--1964), 
\'Etienne Hirsch (1901--1994),
Lew Kowarski  (1907--1979), 
Andr\'e  Labarthe  (1902--1967),
Jean Morin (1897--1943),  
Abraham Robinson (1918--1974),
Yves Rocard (1903--1992),
Alberte and Georges Ungar (1913--2005 and 1906--1977), for example,
one is struck by how wide-ranging were the talents that had managed to escape 
the German occupation of France~\cite{FF}.

Halban and Kowarski brought not only knowledge of nuclear fission, but also
36 gallons of heavy water and a gramme of radium, from Joliot-Curie's laboratory.
An editorial~\cite{N40} in Nature proclaimed that
``French scientific workers who have succeeded in reaching 
Great Britain are eager to help us in our war effort, and we 
must in our turn help them in every possible way." 
Labarthe and Haldane spoke at a three-day 
conference on ``Science and World Order" held in London at the Royal 
Institution, from 26 to 28 September 1941, as did Jacques M\'etadier (1893--1986),
Hirsch (using the alias Bernard) and ``a French man of science -- who desired to remain anonymous''~\cite{SWO}.
Gu\'eron began working at Imperial College London and also
with a Free French sabotage laboratory: {\em la laboratoire de
chimie de l'armement des forces fran\c{c}aises libres},
moving to Cambridge in December 1941, to join Halban and 
Kowarski in the ``Tube Alloys" team working on the atomic bomb~\cite{JG}. 
Rocard was de Gaulle's director of naval research. Like Furry,  he worked on radar.
Morin was director of armaments for the Free French.
Cathala worked on explosives at  the Royal Ordnance Factory.
The Ungars worked in Oxford  with Solly Zuckerman (1904--1993)
on the traumatic shock of bombing. Robinson, a mathematician, enlisted in the Free French Air Force and
was sent in 1941 to the Royal Aircraft Establishment, where he worked on supersonic flow~\cite{AR}.
We have not been able to decide which, if any, of these might have been FFS. In terms of location
and desire for anonymity in October 1941, Gu\'eron is a plausible candidate.

\section{Asymptotics and distributions}

At large $u$, the Dickman function $\rho(u)$ becomes very small.
For example, $\rho(100)\approx1.0006\times10^{-229}$ gives the tiny probability that a
random number of very large size $N$ has no prime divisor $p>N^{1/100}$. After using
Algorithm~1, with $N=101$ steps at $D=350$ decimal digits of working precision, we were well 
equipped to compute more than 100 good digits of $\rho(u)$, for any real $u\in[1,101]$, notwithstanding
extreme cancellations between the polylogarithms of weights $k\leq100$ in the alternating sum
$\rho(u)=\sum_{u>k\ge0}(-1)^kP_k(u)$ over the Furry probabilities $P_k(u)$ that also
solve the Buchstab problem $(u+1)\omega(u+1)=\sigma(u)=\sum_{u>k\ge0}P_k(u)$.

For $u\ge6$, we have $\rho(u)<u^{-u}$ from~\cite{dBr3} and hence $a(u)=\log(\rho(u))/u+\log(u)<0$.
Fig.~\ref{fig:rho} shows the decrease of $a(u)$, with 
$a(100)\approx\log(1.0006\times10^{-29})/100\approx-0.66774$.

\subsection{Sum rule}

Also of interest is the integral $I(u)=\int_0^u\rho(x){\rm d}x$, whose asymptote~\cite{JL} , 
$I(\infty)=e^\gamma$, is the {\em reciprocal} of the Mertens constant in~(\ref{M3}) that 
preoccupied Lindemann and Furry in the Buchstab problem. This leads to the sum rule
\begin{equation}\sum_{n>0}n\rho(n)=\sum_{n>k\ge0}(-1)^knP_k(n)=e^\gamma\label{SR}\end{equation}
obtained by using the Dickman equation~(\ref{DE}) to express an integral of $\rho(x)$ as a sum of values 
$n\rho(n)$ determined by the Furry constants $P_k(n)$ that were stored,
for $k<n\in[1,101]$, in a single run of Algorithm~1. It affords a strong test of the performance
of {\tt polylogmult}, at weights up to 100 and 350-digit precision, giving more than 100 good digits of
the tiny residual integral $\int_{100}^\infty\rho(x){\rm d}x\approx1.5438\times10^{-230}$.

Taking a logarithm, we define $b(u)=\log(e^\gamma-I(u))/u+\log(u)$,
with $0>a(u)>b(u)$ for $u\ge6$ and $b(100)\approx\log(1.5438\times10^{-30})/100\approx-0.68643$.
Fig.~\ref{fig:smr} gives a plot of $a(u)-b(u)$ for $u\in[6,101]$.

\subsection{Oscillations of the Mertens discrepancy}

Following Furry, we study the Mertens discrepancy $\Delta(u)=(u+1)e^{-\gamma}-\sigma(u)$.
Lindemann remarked that the primes make  $\Delta(1)=2e^{-\gamma}-1\approx 0.12292$ positive. 
Furry showed that rough semiprimes make $\Delta(2)=3e^{-\gamma}-1-\log(2)\approx-8.7687\times10^{-3}$ negative.
His polylogarithm $P_2(3)=\log(2)\log(\frac32)-\frac12{\rm Li}_2(\frac14)$, at weight 2 for rough numbers with 3 prime divisors,
gives the small positive value 
\begin{equation}\Delta(3)=4e^{-\gamma}-1-\log(3)-P_2(3)\approx4.9686\times10^{-6}.\label{FD3}\end{equation}
Let $u_n$ be the $n$-th solution to $\Delta(u)=0$ with $u>1$. Then
$u_1\approx1.4833$ solves $(u+1)e^{-\gamma}=1+\log(u)$ and
$u_2\approx2.2270$ solves $(u+1)e^{-\gamma}=1+\log(u)+P_2(u)$,
with a dilogarithm in~(\ref{P2}), The next three zeros
occur at $u_3\approx3.0017$, $u_4\approx3.7858$ and $u_5\approx4.5665$.
Their locations can be found at high precision using Theorems~4 and~5.
Thereafter we use Algorithm~2, which gives 250 good digits, for $u_6\approx5.3507$
up to $u_{25}\approx20.776$, in about a minute. The extrema are easy to locate,
since $u\Delta^\prime(u)=\Delta(u-1)$ for $u>1$. The minimum value of $\Delta(u)$ is
$\Delta(e^\gamma)=e^{-\gamma}-\gamma=-0.015756$. Thereafter, diminishing
local extrema occur at $u=u_n+1$.

The oscillations of $\Delta(u)$ are damped super-exponentially, with
$|\Delta(u)| <u^{-u}$ for $u\ge1$. In Fig.~\ref{fig:osc}
we have divided $\Delta(u)$ by $\rho(u+3)$, to moderate the damping.

\subsection{Filtration by weight}

\begin{table}[h!]\centering
\begin{tabular}{|l|rrrrrrrrrr|}\hline
$u\,\backslash k$ & 0 & 1 & 2 & 3 & 4 & 5 & 6 & 7 & 8 & 9 \\\hline
2 & 591 & 409 &  &  &  &  &  &  &  &  \\
3 & 445 & 489 & 66 &  &  &  &  &  &  &  \\
4 & 356 & 494 & 145 & 5 &  &  &  &  &  &  \\
5 & 297 & 478 & 203 & 22 & 0 &  &  &  &  &  \\
6 & 254 & 456 & 243 & 44 & 2 & 0 &  &  &  &  \\
7 & 223 & 433 & 271 & 67 & 6 & 0 & 0 &  &  &  \\
8 & 198 & 412 & 291 & 88 & 11 & 1 & 0 & 0 &  &  \\
9 & 178 & 391 & 304 & 107 & 18 & 1 & 0 & 0 & 0 &  \\
10 & 162 & 373 & 313 & 125 & 25 & 3 & 0 & 0 & 0 & 0 \\
20 & 85 & 254 & 315 & 218 & 95 & 27 & 5 & 1 & 0 & 0 \\
30 & 57 & 195 & 287 & 246 & 139 & 55 & 16 & 4 & 1 & 0 \\
40 & 43 & 160 & 261 & 253 & 166 & 79 & 28 & 8 & 2 & 0 \\
50 & 35 & 137 & 239 & 253 & 183 & 97 & 40 & 13 & 3 & 1 \\
60 & 29 & 120 & 221 & 249 & 194 & 112 & 50 & 18 & 5 & 1 \\
70 & 25 & 107 & 206 & 245 & 202 & 124 & 59 & 23 & 7 & 2 \\
80 & 22 & 96 & 193 & 239 & 207 & 134 & 68 & 28 & 9 & 3 \\
90 & 20 & 88 & 182 & 234 & 210 & 142 & 76 & 33 & 12 & 4 \\
100 & 18 & 81 & 173 & 228 & 212 & 149 & 82 & 37 & 14 & 4 \\\hline
\end{tabular}
\caption{\large Rounded contributions to $\sigma(u)$ of $P_k(u)$ in parts per thousand.}
\label{table:round}
\end{table}

Following de Bruijn's astute work~\cite{dBr1,dBr2,dBr3}, written almost a decade after 
the Lindemann--Furry letters~\cite{F42,L41a}, subsequent authors~\cite{CG,JF,FGHM,LW,MZW} have 
studied the behaviour of $\rho(u)$ and $\sigma(u)$, with the latter perhaps disguised by de Bruijn's
definition~\cite{dBr1} of $\omega(u)=\sigma^\prime(u)=\sigma(u-1)/u$, for $u>1$.
We recover information that may have been
overlooked: the filtration by weight begun by Furry.

The filtration $\sigma(u)=\sum_{u>k\ge0}P_k(u)$ appeals
to us, as physicists, because $P_k(u)$ is commendably concrete. It estimates
the density, relative to the primes, of large rough numbers of size $N$
that are products of $k+1$ primes, the smallest of which satisfies $p^{u+1}>N$.
Table~\ref{table:round} rounds $10^3P_k(u)/\sigma(u)$, showing
the contribution by weight, in parts per thousand, for some integer values of $u\in[2,100]$.
It is notable that the modal weight is less than 4 for $u\le100$. Table~\ref{table:msd}
indicates how the mean and standard deviation of the distribution by weight
change with $u$.
\begin{table}[h!]\centering
\begin{tabular}{|l|cccccccccc|}\hline
$u$   & 10 & 20 & 30 & 40 & 50 & 60 & 70 & 80 & 90 & 100 \\\hline
mean & 1.4867 & 2.0916 & 2.4660 & 2.7378 & 2.9513 & 3.1272 & 3.2767 & 3.4068 & 3.5218 & 3.6250 \\
s.d.    & 1.0039 & 1.2399 & 1.3724 & 1.4630 & 1.5313 & 1.5857 & 1.6308 & 1.6692 & 1.7025 & 1.7320 \\\hline
\end{tabular}
\caption{\large Approximate means and standard deviations of weight distributions.}
\label{table:msd}
\end{table}

In the tails of the distributions by weight, we encounter $P_k(k+y)$ with $y=O(1)$.
Such terms are tiny for large $k$ and small $y$, being suppressed by the square of a factorial,
with $(k!)^2P_k(k+y)=y^k(1+O(y))$ and $P_k(k+1)<1/(k!)^2$ for all $k>0$.
At $u=101$, the small Mertens discrepancy $|\Delta(101)|\approx1.2931\times10^{-244}$
is much larger than its weight 100 contribution, namely 
$P_{100}(101)\approx4.3068\times10^{-317}$.
The last 4 terms in the sum over weights give a contribution 
$\sum_{k=1}^4P_{101-k}(101)\approx6.0106\times10^{-248}$ that is
significantly smaller than $|\Delta(101)|$. 

\subsection{Inclusion of weights up to 200}
\begin{table}[h!]\centering
\begin{tabular}{|l|rrrrrrrrrrrr|}\hline
$u\,\backslash k$ & 0 & 1 & 2 & 3 & 4 & 5 & 6 & 7 & 8 & 9 & 10 & 11 \\\hline
110 & 16 & 75 & 164 & 223 & 213 & 155 & 89 & 41 & 16 & 5 & 2 & 0 \\
120 & 15 & 70 & 157 & 218 & 214 & 159 & 94 & 46 & 18 & 6 & 2 & 0 \\
130 & 14 & 66 & 150 & 213 & 214 & 164 & 99 & 49 & 21 & 7 & 2 & 1 \\
140 & 13 & 62 & 144 & 208 & 214 & 167 & 104 & 53 & 23 & 8 & 3 & 1 \\
150 & 12 & 59 & 138 & 204 & 213 & 170 & 108 & 57 & 25 & 9 & 3 & 1 \\
160 & 11 & 56 & 133 & 200 & 213 & 173 & 112 & 60 & 27 & 10 & 3 & 1 \\
170 & 10 & 53 & 129 & 196 & 212 & 175 & 116 & 63 & 29 & 11 & 4 & 1 \\
180 & 10 & 51 & 125 & 192 & 211 & 177 & 119 & 66 & 31 & 12 & 4 & 1 \\
190 & 9 & 49 & 121 & 188 & 210 & 179 & 122 & 69 & 33 & 13 & 5 & 1 \\
200 & 9 & 47 & 117 & 185 & 209 & 181 & 125 & 72 & 35 & 14 & 5 & 2 \\\hline
\end{tabular}
\caption{\large Further contributions to $\sigma(u)$ of $P_k(u)$ in parts per thousand.}
\label{table:round2}
\end{table}
\begin{table}[h!]\centering
\begin{tabular}{|l|cccccccccc|}\hline
$u$   & 110 & 120 & 130 & 140 & 150 & 160 & 170 & 180 & 190 & 200 \\\hline
mean & 3.7185 & 3.8040 & 3.8828 & 3.9558 & 4.0239 & 4.0876 & 4.1475 & 4.2040 & 4.2575 & 4.3083 \\
s.d.   & 1.7582 & 1.7820 & 1.8036 & 1.8235 & 1.8418 & 1.8588 & 1.8746 & 1.8895 & 1.9034 & 1.9166 \\\hline
\end{tabular}
\caption{\large Further means and standard deviations of weight distributions.}
\label{table:msd2}
\end{table}
Tables~\ref{table:round2} and~\ref{table:msd2} extend Tables~\ref{table:round} and~\ref{table:msd}
to $u=200$. At $u=201$ the weight 200 Furry probability $P_{200}(201)\approx5.9733\times10^{-751}$ 
makes a very small contribution to
\begin{align}\rho(201)=\sum_{k=0}^{200}(-1)^kP_k(201)&\approx6.7083\times10^{-534}\label{r201}\\
e^\gamma-\sum_{n=1}^{201}n\sum_{k=0}^{n-1}(-1)^kP_k(n)
&\approx9.1972\times10^{-535}\label{i201}\\
202\,e^{-\gamma}-\sum_{k=0}^{200}P_k(201)&\approx7.7106\times10^{-552}\label{d201}\end{align}
for each of which we have obtained more than 400 significant figures.

\section{Statistics of factorization}

Here we  study ranges
of reasonably large integers, $n\in[N_1, N_2]$, and count  primes, semiprimes
and, crucially, integers with precisely three (not necessarily distinct) prime factors, which for convenience
we call triprimes, as suggested by John Conway (1937--2020).
We choose ranges that are relatively narrow, with $N_2\gg(N_2-N_1)\gg1$.

Counting primes is easier than counting semiprimes or triprimes. 
Consider the modest case with $N_2=10^{24}$
and $N_2-N_1=552750053$, which contains $c_0=10^7$ primes,
the smallest of which is $N_1$. One expects about
$(N_2-N_1)/\log(N_2)\approx10002346$ primes from the prime number theorem.
To count the primes,
we first remove all integers with a prime divisor $p<N_2^{1/4}=10^6$.
This is easily done, by crossing out arithmetic progressions in a bitmap. 
There remain 22463197 rough integers. This accords
fairly well with a Mertens estimate $4e^{-\gamma}c_0\approx22458379$.

There are $c_1+c_2=12463197$ composite integers that get through the sieve.
These comprise $c_1$ semiprimes, $n=p_0p_1$,
with $10^6<p_0\le p_1<10^{18}$, and $c_2$ triprimes, $n=p_0p_1p_2$, with
$10^6<p_0\le p_1\le p_2<10^{12}$.
Prime hunters discard these composite numbers, using a probable-primality test,
such as {\tt ispseudoprime} in {\tt Pari/GP}, or an equivalent in dedicated software,
such as {\tt OpenPFGW}~\cite{PFGW}. Any number that fails such a test is certainly composite~\cite{CP}.
Often no attempt is made to count semiprimes, or triprimes.
In this modest range, it is not hard to do so. There are $c_0=10^7$ primes, 
$c_1=10992988$ semiprimes
and $c_2=1470209$ triprimes.

Furry's letter to Nature~\cite{F42} tells us to expect that
$c_1/c_0\to\log(3)\approx1.0986$ and that
$c_2/c_0\to P_2(3)=\log(2)\log(\frac32)-\frac12{\rm Li}_2(\frac14)\approx 0.14722$,
asymptotically, provided that we make very large counts, $c_k$.
Allowing for variations of order $\sqrt{c_k}$, with limited statistics, Furry's
asymptotic predictions,  based on continuous analysis, compare quite well with our 
modest data for numbers of size merely $10^{24}$, in a range with merely $10^7$ primes.

We repeated this process in a case with $N_2=10^{36}$,
where distinguishing a rough semiprime from a rough triprime requires considerably more work.
Here we chose $N_1=N_2-82890279$, which gives $c_0=10^6$ primes, compared
with an expectation of 999966, from the prime number theorem. Sieving out
integers divisible by primes $p<N_2^{1/4}=10^9$, we were left with
2246578 rough numbers. Of these, $c_1=1099722$ are semiprimes 
and $c_2=146856$ are triprimes. Again, the agreement with Furry is
acceptable, allowing for square-root uncertainties and sub-asymptotic effects.

\section{Comments and conclusions}

After formulating a now proven~\cite{Sound} conjecture~(\ref{db}) on 
Dickman polylogarithms and their asymptotic constants~\cite{DB},
one of us (DB) noted that Furry's dilogarithm ${\rm Li}_2(\frac14)$ in~(\ref{Fk1})
for the Mertens discrepancy~\cite{SW} also appears at weight 2 in the Dickman problem.
This led to historical research that appears, much condensed, in Section~3.
The second author (SO) was led to Theorem~2 via the development of a path-integral approach to prime density.
Our joint work synthesizes
best practice in high-energy physics~\cite{VW} with
an efficient double-tail method~\cite{Akh1,Akh2,HC} for multiple
polylogarithms, for which we offer the following summary.
\begin{enumerate}
\item Theorems~1 and~2 solve the Dickman and Buchstab problems
in parallel, with products of rapidly computable multiple
polylogarithms~(\ref{Mjn}) determining the Furry probabilities 
$P_k(u)$ of pure weight $k<u$ in the terminating series~(\ref{both}).
\item A one-off hour-long payment, for constants $P_k(n)$ with integers  $k<n\in[1,101]$,
is rewarded by fast computation of $P_k(u)$ for real $u>k$ in Algorithm 2, giving 
at least 100 significant  figures for the Dickman and Buchstab
functions with $u\in[1,101]$ and also their filtrations by weight.
\item For small $u-k$, Theorem~3 is efficient. For $k<5$,
Theorems~4 and~5 are efficient. For $k<10$,
one-dimensional quadrature based on 
Theorem~6 is efficient. For $k<u\in[10,201]$,
the efficient $N$-step method of Algorithms~1 and~2 was used.
\item The asymptotic appearances of $e^{-\gamma}$ in the Buchstab
problem and $e^\gamma$ in the Dickman problem provide
strong tests of the accuracy and efficiency of the
Akhilesh--Cohen algorithm in the procedure {\tt polylogmult}
of {\tt Pari/GP} at weights up to 200 with 1000-digit working precision.
\item Tables~1 to~4 show distributions by weight, their
means and standard deviations. Figures~1 to~3 show
trends at large $u$ and oscillations about these.
\item Mindful that mathematics is a human activity, 
we have chronicled the roles of 6 interesting authors, of
at least 4 nationalities, two of whom remain anonymous.
\end{enumerate}

\section*{Acknowledgements}

We thank Kevin Acres, Graham Farmelo and Michael St Clair Oakes, for close reading of preliminary drafts,
Steven Charlton and Aur\'elien Dersy, for generous technical advice,
Clare Kavanagh, for valuable help with Section 3,  and Henri Cohen,
for the efficiency of {\tt polylogmult}.
As mathematical physicists who gained their doctorates in England and Germany,
in peaceful circumstances, we are grateful to our fellow physicists,
Lindemann and Furry, for exercising sound, calm
and dispassionate judgement, on a mathematical question of considerable interest,
when our countries were unfortunately at war.

}\raggedright

\newpage

\begin{figure}
\centering
\caption{\large Plot of $a(u)=\log(\rho(u))/u+\log(u)<0$ for $u\in[6,101]$.}
\includegraphics[width=0.9\textwidth]{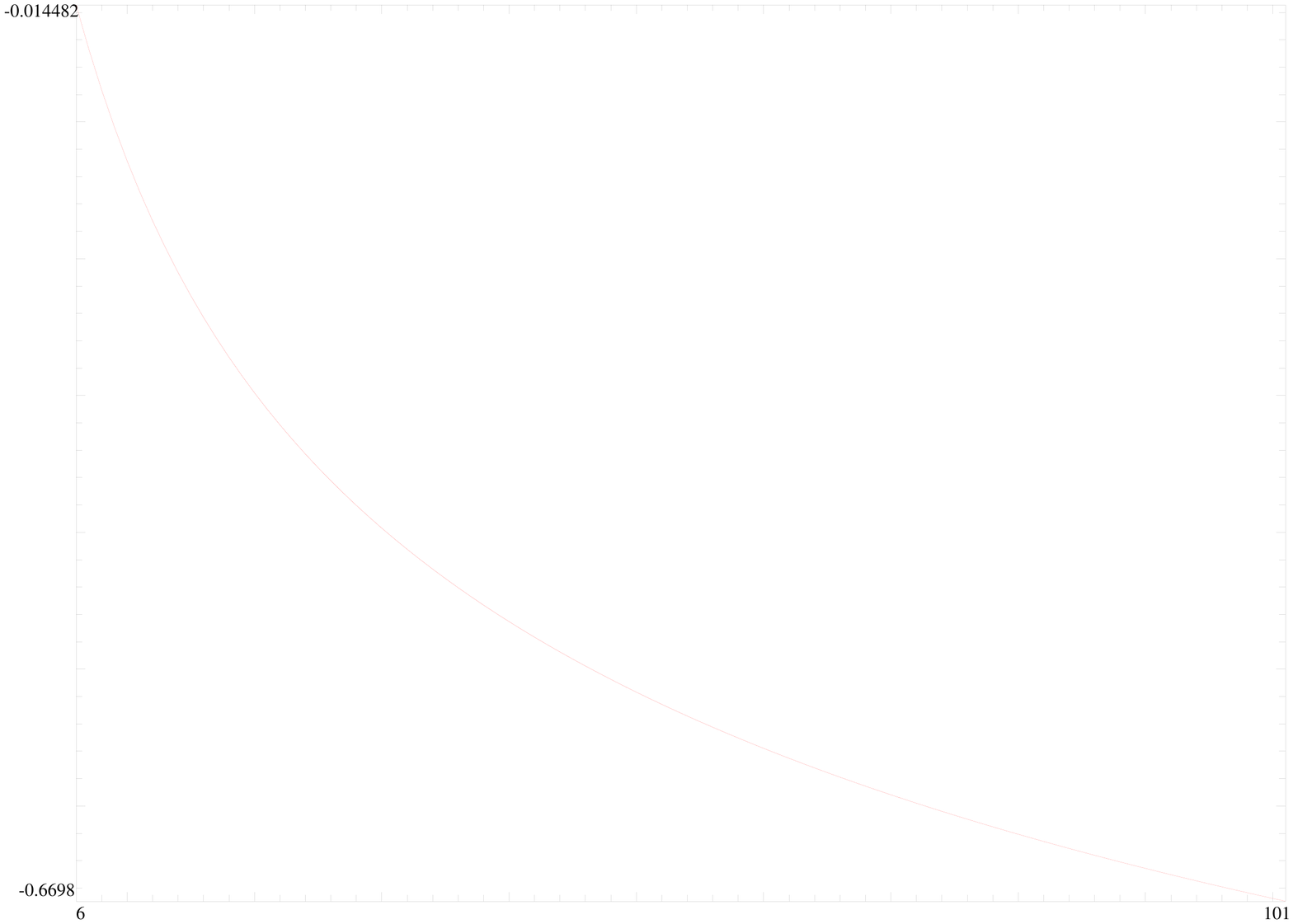}
\label{fig:rho}
\end{figure}

\begin{figure}
\centering
\caption{\large Plot of $a(u)-b(u)=(\log(\rho(u))-\log(\int_u^\infty\rho(x){\rm d}x))/u$ for $u\in[6,101]$.}
\includegraphics[width=0.9\textwidth]{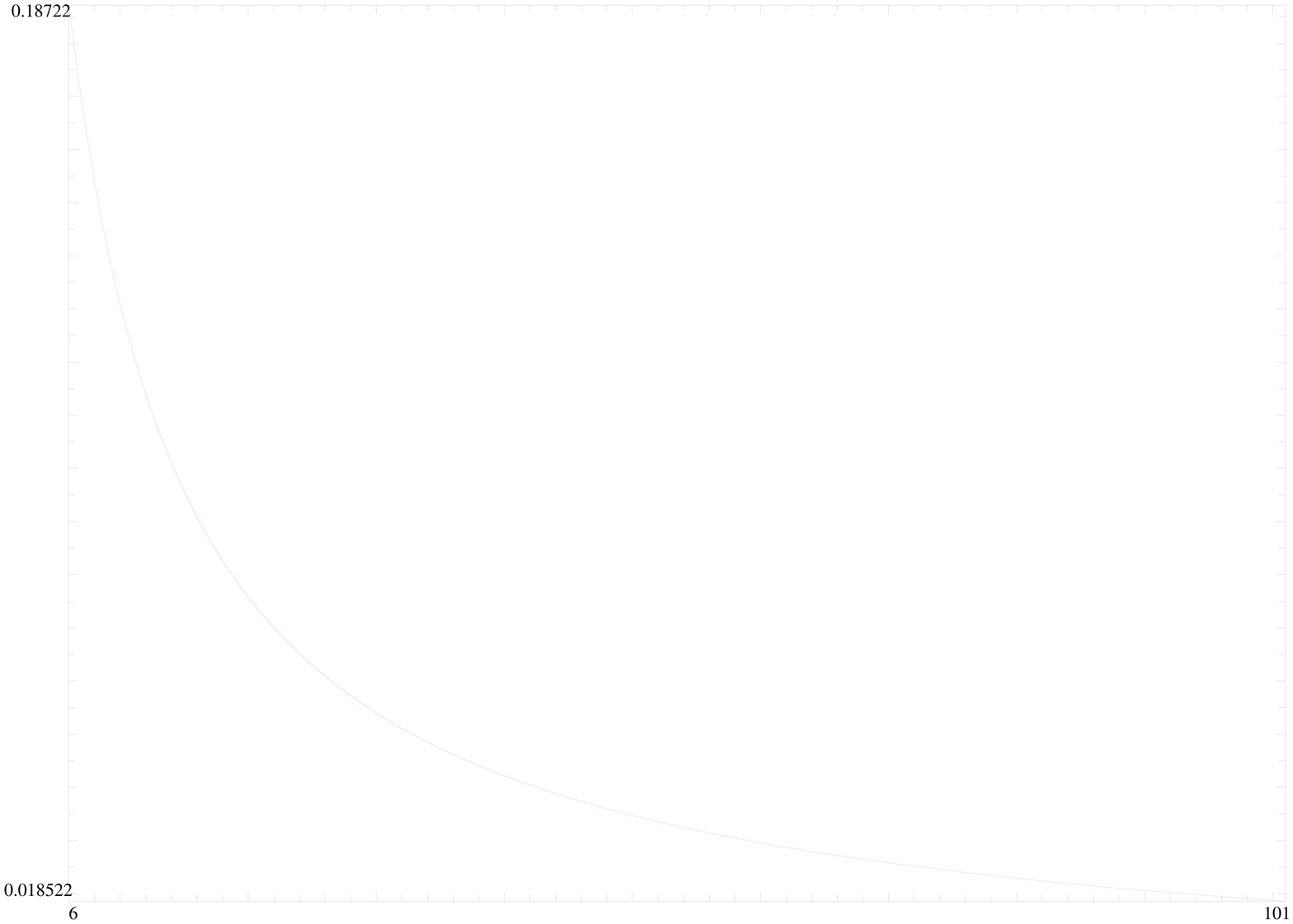}
\label{fig:smr}
\end{figure}

\begin{figure}
\centering
\caption{\large Oscillations of $\Delta(u)/\rho(u+3)$ for $u\in[6,21]$.}
\includegraphics[width=0.9\textwidth]{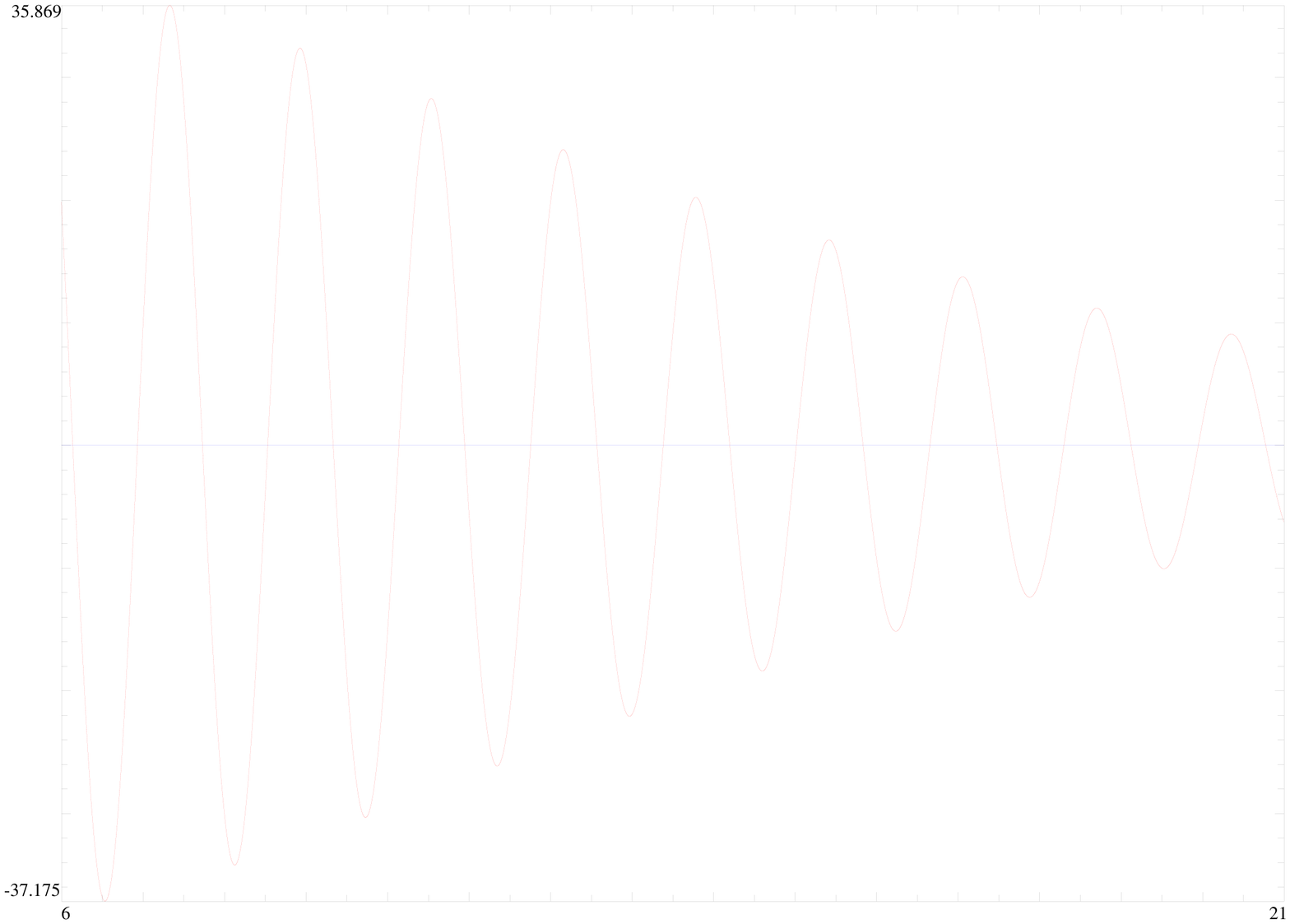}
\label{fig:osc}
\end{figure}

\end{document}